\documentclass[11pt]{article}
\usepackage{amsfonts,amsmath,tikz}
\usepackage{epsfig}
\usepackage{latexsym}
\usepackage{color}
\usepackage{epsf,amssymb}

\usepackage{amsthm}
\usepackage{yfonts}
\usepackage{hyperref}
\usepackage{ifthen}
\usepackage{enumitem,enumerate}
\usepackage[T1]{fontenc} 
\usepackage{diagbox}
\usepackage{array}
\usetikzlibrary{arrows.meta}
\usepackage[normalem]{ulem}
\usepackage{tabularx}

\usepackage[ruled,vlined]{algorithm2e}
\SetKw{And}{and}
\SetKw{Or}{or}

\usepackage{float} 
\floatstyle{plaintop}
\restylefloat{table} 

\newtheorem{theorem}{Theorem}
\newtheorem{lemma}[theorem]{Lemma}

\newtheorem{corollary}[theorem]{Corollary}
\newtheorem{proposition}[theorem]{Proposition}

\newtheorem{question}{Question}
\newtheorem{remark}{Remark}

\def\vertex(#1){\put(#1){\circle*{2}}}
\def\vertexo(#1){\put(#1){\circle{2}}}
\def\vert(#1){\put(#1){\circle*{1.5}}}
\def\verto(#1){\put(#1){\circle{1.5}}}
\def\lab(#1)#2{\put(#1){\makebox(0,0)[c]{#2}}}
\newcommand{\es}{\alpha_{\rm es}}
\newcommand{\al}{\alpha}

\tikzset{My Style/.style={draw, circle, fill=black,scale=0.3}} 

\title{Externally supported independence in split and proper interval graphs}

\author{Dragana Božović\textsuperscript{a,b} \and 
Tanja Dravec\textsuperscript{b,c} \and
Aleksander Kelenc\textsuperscript{a,b}
}     

\begin{document}
\maketitle

\begin{center}
\textsuperscript{a}University of Maribor, FERI, Koroška cesta 46, 2000 Maribor, Slovenia \\
\textsuperscript{b}University of Maribor, FNM, Koroška cesta 160, 2000 Maribor, Slovenia \\
\textsuperscript{c}IMFM, Jadranska ulica 19, 1000 Ljubljana, Slovenia \\
\end{center}

\begin{abstract}
An independent set $B$ in a graph $G=(V,E)$ satisfying the additional condition that the vertices in $N(B)$ are dominated by vertices in $V(G) \setminus N[B]$ is called an externally supported independent set (ESI-set for short). The maximum cardinality of such a set in $G$ is called externally supported independence number of $G$. We study the externally supported independence number in split graphs and proper interval graphs, two graph classes with unbounded clique-width. For split graphs we investigate the structure and computational complexity of the \textsc{ESI Set} problem. We prove that the \textsc{ESI Set} problem remains NP-complete when restricted to split graphs, while also establishing simple upper and lower bounds for the externally supported independence number. These bounds are sharp for a wide range of split graphs.
For split graphs, we further investigate the relationship between the externally supported independence number and the isolation number. In particular, we characterize when an optimal ESI-set yields an optimal isolating set, and conversely, when an optimal isolating set yields an optimal ESI-set. Finally, we consider proper interval graphs and present a linear-time algorithm for finding a maximum ESI-set of an arbitrary proper interval graph.
\end{abstract}

\noindent
{\bf Keywords:} externally supported independent set, externally supported independence number, isolation number, split graphs, proper interval graphs

\medskip



\section{Introduction}

Domination is one of the central topics in graph theory, with numerous variants and generalizations arising from the basic question of selecting a small set of vertices that has a prescribed influence on the remaining vertices; see, for example, \cite{HaynesHedetniemiHenning2022}. One such variant is \emph{isolation}. An early form of this concept was introduced by Lewis et al.~\cite{LewisEtAl2010} under the name \emph{vertex-edge domination}. The isolation number in its present form was introduced and systematically studied by Caro and Hansberg~\cite{CaroHansberg2017}. Instead of requiring every vertex to be dominated, an isolating set is allowed to leave some vertices undominated, provided that they induce no edges. More precisely, a set $A\subseteq V(G)$ is an isolating set if $G\setminus N[A]$ is edgeless, and the minimum cardinality of such a set is the \emph{isolation number}, denoted by $\iota(G)$. Since its introduction, isolation has attracted considerable attention. General bounds and extremal questions have been studied in \cite{CaroHansberg2017,Zylinski2019}, while particular attention has been given to special graph classes, including cubic graphs \cite{ZiemannZylinski2020}, maximal outerplanar graphs \cite{TokunagaEtAl2019}, unicyclic and block graphs \cite{LemanskaEtAl2024}, and graphs attaining the general upper bound on the isolation number \cite{BoyerGoddard2024}. More recently, the behavior of the isolation number under graph products and on generalized Sierpiński graphs has also been investigated \cite{brevsar2027isolation}. Several extensions and variants of isolation have been introduced as well, making isolation a rapidly developing branch of domination theory.

A useful way of viewing an isolating set is to focus not on the vertices selected into the set, but rather on the vertices that remain after its closed neighborhood has been removed. If $A$ is an isolating set of a graph $G$, then $B=V(G)\setminus N[A]$ is an independent set. This observation suggests reversing the usual perspective: instead of starting with a set that isolates the graph, one may ask which independent sets can occur as the vertices left undominated by an isolating set. This point of view led to the introduction of \emph{externally supported independent sets} in \cite{1}. More precisely, an independent set $B\subseteq V(G)$ is an \emph{externally supported independent set}, or an \emph{ESI-set}, if every vertex in $N(B)$ has a neighbor in $D=V(G)\setminus N[B]$.
The maximum cardinality of an ESI-set is called the \emph{externally supported independence number} of $G$ and is denoted by $\alpha_{es}(G)$ \cite{1}. Thus, an ESI-set is not merely an independent set: every vertex adjacent to the set must have a ``support'' outside its closed neighborhood. Equivalently, $N(B)$ must be dominated from the part of the graph lying outside $N[B]$. This additional condition reflects precisely the structural restriction imposed on an independent set when it occurs as the remainder of an isolating set. The parameter $\alpha_{es}(G)$ provides a natural link between independence and isolation. In particular, if $A$ is an isolating set and $B=V(G)\setminus N[A]$, then $B$ is an ESI-set. Consequently, the externally supported independence number yields the lower bound
$
\iota(G)\geq
\left\lceil
\frac{n(G)-\alpha_{es}(G)}{\Delta(G)+1}
\right\rceil,
$
which may improve the corresponding bound formulated in terms of the ordinary independence number \cite{1}. The initial study of $\alpha_{es}(G)$ established several sharp bounds, investigated its relation to the independence number, characterized graphs with $\alpha_{es}(G)=0$, and obtained refined results for trees. From the computational point of view, the \textsc{ESI Set} problem was shown to be NP-complete for general graphs, whereas a linear-time algorithm was obtained for trees \cite{1}.

These results naturally raise the question of identifying graph classes on which the externally supported independence number can be computed efficiently. Indeed, this question was explicitly posed in \cite{1}. An important general observation is that the \textsc{ESI Set} problem can be expressed in monadic second-order logic, see~\cite{1}. Consequently, Courcelle's theorem implies linear-time solvability on every graph class of bounded treewidth \cite{Courcelle1990,Courcelle1992}. Moreover, the Courcelle--Makowsky--Rotics meta-theorem implies polynomial-time solvability of the corresponding optimization problem on graph classes of bounded clique-width \cite{CourcelleMakowskyRotics2000}. In particular, since trees have treewidth one, the existence of a linear-time algorithm for trees also follows from the general MSO framework, although an explicit dynamic programming algorithm was given in \cite{1}.
A further question is what happens on graph classes of unbounded clique-width. For such classes, the general meta-theorems above no longer provide an automatic existence of polynomial-time algorithm, and the complexity of a particular MSO-definable problem has to be investigated using more specific structural properties of the graphs under consideration. In this paper, we study this question for two important subclasses of chordal graphs, namely split graphs and interval graphs, with particular emphasis on proper interval graphs.

Split graphs were introduced by Földes and Hammer~\cite{FoldesHammer1977} and can be characterized in several equivalent ways; in particular, they are precisely the $(2K_2,C_4,C_5)$-free graphs. They form a classical hereditary graph class with a particularly simple structure: their vertex set can be partitioned into a clique and an independent set. This decomposition makes split graphs a natural setting for studying combinatorial and algorithmic problems, since many structural properties can be expressed directly in terms of the interaction between the two parts. At the same time, split graphs have unbounded clique-width~\cite{MakowskyRotics1999}; see also~\cite{BrandstadtEtAl2016}. Hence, the tractability of the \textsc{ESI Set} problem on split graphs does not follow from the bounded-clique-width framework. In fact, we show that the \textsc{ESI Set} problem is NP-complete on split graphs. Since split graphs are chordal, this also establishes NP-completeness for the class of chordal graphs.

We next consider interval graphs, another fundamental subclass of chordal graphs. An interval graph is the intersection graph of a family of intervals on the real line: its vertices correspond to intervals, and two vertices are adjacent precisely when the corresponding intervals intersect~\cite{Brandstadt1999,McKee1999}. Interval graphs have been extensively studied from structural and algorithmic perspectives and have found applications in areas such as scheduling, planning, and computational biology~\cite{Brandstadt1999,McKee1999,Waterman1986,Coulon2006}. Despite their strong structural properties, interval graphs have unbounded clique-width~\cite{GolumbicRotics2000}. They therefore provide a natural class in which to investigate whether the specific geometry of an interval representation can yield efficient algorithms for the \textsc{ESI Set} problem beyond what is guaranteed by the general clique-width meta-theorem.
A particularly structured subclass of interval graphs is given by proper interval graphs. A graph is a proper interval graph if it admits an interval representation in which no interval properly contains another; equivalently, proper interval graphs are precisely the unit interval graphs~\cite{roberts1969,gardi2007}. Proper interval graphs have been extensively studied and possess a rich structural theory as well as efficient recognition and representation algorithms. Interestingly, despite their additional structural restrictions, proper interval graphs also have unbounded clique-width~\cite{GolumbicRotics2000,lozin2008}. Thus, they constitute another natural setting in which the complexity of an MSO-definable problem cannot be inferred directly from bounded clique-width. In this setting, however, the additional structure of proper interval representations can be exploited.

In the next section we present some terminology and preliminary results. Then in Section~\ref{s:split} we initiate a systematic study of the externally supported independence number on split graphs. We first investigate the structure of ESI-sets with respect to a split partition and show that ESI-sets in split graphs fall into two distinct types. ESI-sets of the first type yield a straightforward bound on the externally supported independence number of split graphs. In contrast, ESI-sets of the second type are substantially more difficult to detect; in fact, their existence is closely related to the computational hardness of the problem and is used to prove that the \textsc{ESI Set} problem remains NP-complete even when restricted to split graphs. We further exploit the structural description of ESI-sets to characterize split graphs $G$ attaining the extremal values $\es(G)=\alpha(G)-1$ and $\es(G)=0$, corresponding to the general upper and lower bounds for the externally supported independence number. Finally, we use this structure to investigate more closely the relationship between externally supported independence number and isolation number of a split graph. In particular, we characterize split graphs for which there exists an $\es(G)$-set $B$ such that $V(G)\setminus N[B]$ is an $\iota(G)$-set, as well as split graphs for which there exists an $\iota(G)$-set $D$ such that $V(G)\setminus N[D]$ is an $\es(G)$-set. In Section~\ref{s:interval} we present a linear-time algorithm that returns an $\es(G)$-set of an arbitrary proper interval graph.


\section{Preliminaries}
Let $G$ be a graph on $n(G)$ vertices. 
For a vertex $v\in V(G)$, let $N_G(v)$ denote the \emph{open neighborhood} of $v$, and define its \emph{closed neighborhood} by $N_G[v]=N_G(v)\cup \{v\}$. For any subset $A\subseteq V(G)$, we write $N_G[A]=\bigcup_{a\in A}N_G[a]$ and $N_G(A)=N_G[A]\setminus A$. Whenever the underlying graph is clear from the context, the subscript $G$ will be omitted, and we will simply use the notation $N(v)$, $N[v]$, $N(A)$, and $N[A]$. For a set $S\subseteq V(G)$ and $x \in V(G)$, $N_S(x)$ denotes the set of neighbors of $x$ that are contained in $S$, i.e.\ $N_S(x)=\{y \in S:\, xy \in E(G)\}$. A set $S \subseteq V(G)$ is called a \emph{clique} if any two vertices of $S$ are adjacent in $G$. The \emph{clique number} of $G$, $\omega(G)$, is the cardinality of a largest clique in $G$. A set $P\subseteq V(G)$ is a $k$-\emph{packing} of $G$ if the distance between any two distinct vertices of $P$ is more than $k$, and the $k$-\emph{packing number} of $G$, $\rho_k(G)$, is the maximum cardinality of a $k$-packing of $G$. Clearly, $\rho_1(G)=\alpha(G)$. A set $D \subseteq V(G)$ is a \emph{dominating set} if $N_G[D]=V(G)$ and is an \emph{isolating set} if $G\setminus N_G[D]$ is an empty graph. The cardinality of a smallest dominating set is the \emph{domination number} of $G$, $\gamma(G)$, and the cardinality of a minimum isolating set is called  the \emph{isolation number} of $G$ and it is denoted by $\iota(G)$. 

Recently, the externally supported  independence number was introduced, as a graph parameter that improves the well-known lower bound for the isolation number of a graph~\cite{1}. A set $B \subseteq V(G)$ is called an \emph{externally supported independent set} (an \emph{ESI--set} for short) if $B$ is an independent set and every $v\in N(B)$ has a neighbor in $D=V(G) \setminus N[B]$. The \emph{externally supported independence number} of $G$ (shortly the \emph{ESI--number}), denoted by $\es(G)$, is the maximum cardinality of an ESI-set in $G$. An ESI-set of cardinality $\es(G)$ will be referred to as an $\es(G)$-set. For convenience, we shall refer to the requirement that every vertex in $N(B)$ has a neighbor in $D$ as the \emph{ES-condition}. The definition directly implies the following.
\begin{remark}\label{r:basic}
If $B$ is an ESI-set of $G$ and $D=V(G)\setminus N[B]$, then $D$ is an isolating set of $G$ and if $D$ is an isolating set of $G$, then $V(G)\setminus N[D]$ is an ESI-set of $G$. 
\end{remark}

On the other hand, the ESI-number of $G$ cannot necessarily be obtained from the isolation number of $G$, as for an $\iota(G)$-set $D$, the set $V(G) \setminus N[D]$ is not necessarily the largest ESI-set of $G$. In this paper we will also consider when optimal ESI-sets can be obtained from optimal isolating sets. For this reason we start with the following trivial result.

\begin{proposition}\label{optimalToOptimal}
Let $G$ be a graph for which there exists an $\es(G)$-set $B$ such that $V(G)\setminus N[B]$ is an $\iota(G)$-set. Then there exists $\iota(G)$-set $D$ such that $V(G)\setminus N[D]$ is an $\es(G)$-set.
\end{proposition}
\begin{proof}
Let $B$ be an $\es(G)$-set such that $D=V(G)\setminus N[B]$ is an $\iota(G)$-set. By the ES-condition every vertex in $N(B)$ has at least one neighbor in $D$ and thus $N[D]=D \cup N(B)$. Hence $D$ is an $\iota(G)$-set such that $V(G) \setminus N[D]=B$ is an optimal ESI-set.
\end{proof}

Next, we present some results proved in~\cite{1} that will be needed in the rest of the paper.

\begin{theorem}\cite{1} 
\label{alpha_up_bound}
For any graph $G$ with at least one edge we have 
$$0\leq\es(G)\leq\al(G)-1.$$

\end{theorem}

\begin{lemma}\cite{1}
\label{notinESI}
Let $G$ be a graph and let $x\in V(G)$. Then  $x$ does not belong to any ESI-set of $G$ if and only if there exists $y\in N(x)$ such that $N[y]\subseteq N[x]$.
\end{lemma}

In~\cite{1} ESI-number was computed for some standard graph families. It was also observed that $\es(G) \geq \rho_3(G)$ holds for any graph $G$ with no components of order at least 3. Moreover, the problem of describing all graphs having $\es(G)=\rho_3(G)$ was also posed in~\cite{1}. The following straightforward result on packings implies that paths and cycles are contained in this class and it also gives the ESI-number of paths and cycles.

\begin{proposition}\label{p:k_independence}
If $G$ is a graph of maximum degree 2 having no component isomorphic to a complete graph, then $\es(G)=\rho_3(G)$.
\end{proposition}
\begin{proof}
We already know that $\es(G)  \geq \rho_3(G)$. To prove the reversed inequality, let $B$ be an $\es(G)$-set of a graph $G$ with $\Delta(G)=2$. Since $B$ is an ESI-set, it is clearly an independent set. Hence $d_G(b_1,b_2) \geq 2$ for any two distinct vertices $b_1,b_2 \in B$. Suppose that there exist $b_1,b_2 \in B$ with $d_G(b_1,b_2) \leq 3$. Let $x$ be the neighbor of $b_1$ in a shortest $b_1,b_2$-path. Then $x \in N(B)$ and has no neighbors in $G\setminus N[B]$, as $\Delta(G)=2$. This contradicts the fact that $B$ is an ESI-set of $G$. Hence $d_G(b_1,b_2) \geq 4$ for any two distinct vertices $b_1,b_2 \in B$ and thus $B$ is a $3$-packing of $G$, i.e.\ $\rho_3(G) \geq \es(G)$.
\end{proof}
As a direct consequence of Proposition~\ref{p:k_independence}
we get that $\es(P_n)=\rho_3(P_n)=\lceil \frac{n}{4} \rceil$ and $\es(C_n)=\rho_3(C_n)=\lfloor \frac{n}{4} \rfloor$.

\section{Split graphs}\label{s:split}

Recall that a graph $G$ is a \emph{split graph} if its vertex set can be partitioned into a clique and an independent set. For a split graph $G$, we write $G=K+I$ whenever the vertex set of $G$ admits a partition into a clique $K$ and an independent set $I$ such that $|I|=\alpha(G)$.

We will first present the structure of ESI-sets in split graphs.

\begin{lemma}\label{l:insideI}
Let $G=K+I$ be a split graph with $|I| = \alpha(G)$ and let $B$ be an ESI-set of $G$. Then $B \subseteq I$.
\end{lemma}
\begin{proof}
Suppose that $B \cap K \neq \emptyset$. Hence there exists $x \in B \cap K$. Since $|I|=\alpha(G)$, every vertex in $V(G) \setminus I=K$ has at least one neighbor in $I$, thus there exists $x' \in I$ such that $xx' \in E(G)$. 
Hence $N[x'] \subseteq N[x]$ and thus by Lemma~\ref{notinESI}, $x$ does not belong to any ESI-set of $G$, a contradiction.

\end{proof}

\begin{lemma}\label{l:type1ESIsets}
If $G=K+I$ is a split graph with $|I| = \alpha(G)$, then for any $x \in K$, the set $B=I \setminus N_I(x)$ is an ESI-set of $G$.
\end{lemma}
\begin{proof}
 To show that $B$ is an ESI-set we have to prove the ES-condition, as $B$ is clearly an independent set. Thus let $y$ be an arbitrary vertex in $N(B) \subseteq K$. Hence $xy \in E(G)$ and since $x \notin N(B)$, $x$ is a neighbor of $y$ contained in $V(G)\setminus N[B]$, or equivalently $y$ satisfies the ES-condition. 
\end{proof}

In a split graph $G=K+I$ with $|I|=\alpha(G)$ an ESI-set $B=I \setminus N_I(x)$ for some $x \in K$ described in Lemma~\ref{l:type1ESIsets} is called a \emph{type I} ESI-set.

\begin{corollary}
If $G=K+I$ is a split graph with $|I| = \alpha(G)$, then $$\es(G) \geq n(G)-\min{\{{\rm{deg}}(x):\, x \in K\}}-1.$$
\end{corollary}

\begin{proof}
Let $G=K+I$ be a split graph with $|I| = \alpha(G)$ and let $x \in K$ be a vertex with a smallest degree, i.e.\ deg$(x)=\min{\{{\rm{deg}}(y):\, y \in K\}}$. By Lemma~\ref{l:type1ESIsets} $B=I \setminus N_I(x)$ is an ESI-set. Since $|B|=n(G)-\min{\{{\rm{deg}}(x):\, x \in K\}}-1$, the bound follows. 
\end{proof}

\begin{lemma}\label{l:type2ESIsets}
Let $G=K+I$ be a split graph with $|I| = \alpha(G)$ and let $B \subseteq I$ with $N(B)=K$. Then $B$ is an ESI-set of $G$ if and only if $N(I\setminus B)=K$. 
\end{lemma}

\begin{proof}
Assume first that $B \subseteq I$ is an ESI-set with $N(B)=K$. By ES-condition any $x \in N(B)=K$ has a neighbor in $V(G) \setminus N[B]=I\setminus B$. Thus $N(I \setminus B)=K$ as desired. Conversely let $B \subseteq I$ be a set with $N(B)=K$ such that $N(I \setminus B)=K$. Hence, any $x \in K=N(B)$ has a neighbor in $I \setminus B=V(G) \setminus N[B]$ and thus $B$ satisfies the ES-condition. Since $B$ is also independent, it is an ESI-set of $G$. 
\end{proof}

In a split graph $G=K+I$ with $|I|=\alpha(G)$ an ESI-set $B \subseteq I$ with $N(B)=N(I \setminus B)=K$ described in Lemma~\ref{l:type2ESIsets} is called a \emph{type II} ESI-set.

In the next result we present a formula for computing the ESI-number of an arbitrary split graph.

\begin{theorem}\label{thm:split}
If $G=K+I$ is a split graph with $|I|=\alpha(G)$ that is not isomorphic to a complete graph, then $$\es(G)=\max{\{|I'|:\, I' \subseteq I, N(I') \subset K \textrm{ or } N(I') = N(I \setminus I')=K\}}.$$
\end{theorem}
\begin{proof}
Let $B$ be an $\es(G)$-set. By Lemma~\ref{l:insideI}, $B \subseteq I$. Suppose that $$|B| > \max{\{|I'|:\, I' \subseteq I, N(I') \subset K \textrm{ or } N(I') = N(I\setminus I')=K\}}.$$ Hence $N(B)=K$ and there exists $x \in N(B)\setminus N(I \setminus B)$. Then $x \in N(B)$ and it is not dominated by $G\setminus N[B]$, a contradiction with $B$ being an ESI-set. Hence $|B| \leq \max{\{|I'|:\, I' \subseteq I, N(I') \subset K \textrm{ or } N(I') = N(I\setminus I')=K\}}$. By Lemma~\ref{l:type1ESIsets} it follows that the largest possible type I ESI-set of $G$ has cardinality $\max{\{|I'|:\, I' \subseteq I, N(I') \subset K\}}$ and the largest type II ESI-set has cardinality $\max{\{|I'|:\, I' \subseteq I, N(I') = N(I\setminus I')=K\}}$ by Lemma~\ref{l:type2ESIsets}. This  implies the existence of an ESI-set of cardinality 
$$\max{\{|I'|:\, I' \subseteq I, N(I') \subset K \textrm{ or } N(I')= N(I\setminus I')=K\}}.$$
\end{proof}

In~\cite{1} trees that achieve the upper bound from Theorem~\ref{alpha_up_bound} were characterized. From Theorem~\ref{thm:split} we can easily get a characterization of split graphs $G$ with $\es(G)=\alpha(G)-1$.

\begin{corollary}\label{c:split_uper}
If $G=K+I$ is a split graph with $|I|=\alpha(G)$, then $\es(G)=\alpha(G)-1$ if and only if $\omega(G)=|K|+1$ or there exists $y \in K$ with ${\rm{deg}}_I(y)=1$.
\end{corollary}

\begin{proof}
First, let $\es(G)=\alpha(G)-1$ and let $B$ be an $\es(G)$-set of a split graph $G=K+I$ with $|I|=\alpha(G)$. It follows from Lemma~\ref{l:insideI} that $B \subseteq I$. Since $|B|=\alpha(G)-1$, there exists $x \in I$ such that $B=I \setminus \{x\}$. If $N(B)=K$, then $V(G) \setminus N[B]=\{x\}$ and by ES-condition $x$ dominates $N(B)=K$. Hence $N(x)=K$ and thus $\omega(G)=|K|+1$. If $N(B) \subset K$, then there exists $y \in K$ that has no neighbor in $B$ and since $I$ is an $\alpha(G)$-set, $yx \in E(G)$. Thus ${\rm{deg}}_I(y)=1$.

For the converse assume first that $\omega(G)=|K|+1$, i.e.\ there exists $x \in I$ with $N(x)=K$. Then $B=I \setminus \{x\}$ is an ESI-set of $G$ of cardinality $|I|-1=\alpha(G)-1$. Indeed, ES-condition is satisfied, as any vertex in $N(B)$ is dominated by $x \notin N[B]$. Assume now that there exists $y \in K$ with ${\rm{deg}}_I(y)=1$. In this case let $B=I \setminus N_I(y)$. Then it follows from Lemma~\ref{l:type1ESIsets} that $B$ is an ESI-set of $G$. Since $|B|=|I|-1=\alpha(G)-1$, Theorem~\ref{alpha_up_bound} implies that $\es(G)=\alpha(G)-1$.
\end{proof}

Note that there is no nontrivial tree $T$ with $\es(T)=0$ and Theorem~\ref{thm:split}  implies that complete graphs are the only split graphs $G$ with $\es(G)=0$.

\begin{corollary}\label{p:split_lower}
If $G=K+I$ is a split graph not isomorphic to a complete graph, then $\es(G) \geq 1$. 
\end{corollary}
\begin{proof}
Let $x$ be an arbitrary vertex of $I$. If $N(x) \subset K$, then there exists $y \in K \setminus N(x)$. Since $y$ dominates $N(x)$, $\{x\}$ is an ESI-set of $G$ and thus $\es(G) \geq 1$. If $N(x)=K$, then $I \setminus \{x\}$ is an ESI-set which implies that $\es(G) \geq |I|-1 = \alpha(G)-1 \geq 1$.
\end{proof}

Although all ESI-sets in split graphs can be easily described, finding a maximum ESI-set is computationally hard due to ESI-sets of type II.
We next show that the \textsc{ESI Set} problem is NP-complete even when restricted to split graphs.

For the NP-hardness proof, we give a polynomial-time reduction from the well-known \textsc{Set Cover} problem. We first recall the two decision problems.

\begin{center}
\fbox{\parbox{0.9\linewidth}{\noindent
\textbf{\textsc{Set Cover}}\\[.8ex]
\begin{tabular*}{.93\textwidth}{rl}
{\em Instance:} & A finite set $U$, a collection $\mathcal{S}=\{S_1,\ldots,S_m\}$ of subsets of $U$ such \\ & that $\bigcup_{S_i\in\mathcal{S}} S_i=U$, and a positive integer $k$.\\
{\em Question:} & Does there exist a subcollection $\mathcal{C}\subseteq\mathcal{S}$ with $|\mathcal{C}|\leq k$ such that \\ & $\bigcup_{S_i\in\mathcal{C}}S_i=U$?
\end{tabular*}
}}
\end{center}

\begin{center}
\fbox{\parbox{0.9\linewidth}{\noindent
\textbf{\textsc{ESI Set}}\\[.8ex]
\begin{tabular*}{.93\textwidth}{rl}
{\em Instance:} & A split graph $G$ and a positive integer $\ell$.\\
{\em Question:} & Does $G$ contain an ESI-set of size at least $\ell$?
\end{tabular*}
}}
\end{center}

Let $(U,\mathcal{S},k)$ be an instance of \textsc{Set Cover}, where $\mathcal{S}=\{S_1,\ldots,S_m\}$. From this instance, we construct a split graph $G$ as follows. The vertices corresponding to the elements of $U$ induce a clique, while each set $S_i\in\mathcal{S}$ is represented by a vertex $s_i$ in the independent set. A vertex $s_i$ is adjacent to a vertex $u\in U$ if and only if $u\in S_i$.

We construct $G'$ from $G$ by adding $k$ copies $s_i^1,\ldots,s_i^k$ of each vertex $s_i$.
For convenience, let $s_i^0=s_i$ and for any $i \in \{1,\ldots , m\}$ let
$$X_i=\{s_i^0,s_i^1,\ldots,s_i^k\}.$$
All vertices in $X_i$ have the same neighborhood in $U$, that is,
$$N(s_i^j)=S_i\qquad\text{for every }j\in{0,\ldots,k}.$$
Thus, $G'$ is a split graph with clique $U$ and independent set $X=\bigcup_{i=1}^m X_i$.

\begin{proposition}\label{prop:esi-split}
Let $(U,\mathcal{S},k)$ be an instance of \textsc{Set Cover}, where $|\mathcal{S}|=m$. Then $(U,\mathcal{S},k)$ has a set cover of size at most $k$ if and only if $G'$ has an ESI-set of size at least $m(k+1)-k$.
\end{proposition}

\begin{proof}
Suppose that $(U,\mathcal{S},k)$ has a set cover $\mathcal{C}\subseteq\mathcal{S}$ with $|\mathcal{C}|\leq k$.
Define $B=X\setminus \{s_i^0; \, S_i \in \mathcal{C} \}$ which implies that $X\setminus B = \bigcup_{S_i \in {\cal{C}}}\{s_i^0\}$. Then $|B|=m(k+1)-|\mathcal{C}|\geq m(k+1)-k$. Since $B\subseteq X$ and $X$ is an independent set in $G'$, it follows that $B$ is also an independent set. Since $\cal{C}$ is a set cover, it follows that $$N(X \setminus B)=N\left(\bigcup_{S_i \in {\cal{C}}}\{s_i^0\}\right)=U=N\left(\bigcup_{S_i \in {\cal{C}}}\{s_i^1\}\right) = N(B).$$ 
Thus, by Lemma~\ref{l:type2ESIsets}, $B$ is a (type II) ESI-set of $G'$.


Now, let $B$ be an ESI-set of $G'$ of cardinality $|B|\geq m(k+1)-k$. Since $G'$ is obtained from an instance of a set cover, it follows that every vertex from $U$ has at least one neighbor in $X$. Thus, $|X|=\alpha(G')$ and hence by Lemma~\ref{l:insideI}, $B \subseteq X$.
Moreover, since $|X|=m(k+1)$ and $|B|\geq m(k+1)-k$, we have $|X\setminus B|=|X|-|B|\leq k$. Define $\mathcal{C}=\{S_i\in \mathcal{S}; \, X_i\setminus B\neq \emptyset \}$. Then $|\mathcal{C}|\leq |X\setminus B|\leq k$. It remains to prove that $\mathcal{C}$ covers $U$. We first claim that $U\subseteq N(B)$. Since $|X\setminus B|\leq k$ and $|X_i|=k+1$, we have $X_i\cap B\neq \emptyset$ for every $i$. Since $\bigcup_{S_i\in\mathcal{S}}S_i=U$, every $u\in U$ belongs to some $S_i$, and hence has a neighbor in $X_i\cap B$. Therefore, $U\subseteq N(B)$. Let $u\in U$. Since $u\in N(B)$ and $B$ is an ESI-set in $G'$, there exists a neighbor of $u$ in $V(G')\setminus N[B]$. Since $U\subseteq N[B]$, this neighbor must be some $s_i^j\in X_i$. Hence, $s_i^j\notin B$, and therefore $X_i\setminus B\neq\emptyset$, which implies that $S_i\in\mathcal{C}$. Moreover, $us_i^j\in E(G')$ implies that $u\in S_i$. Thus, every $u\in U$ belongs to some $S_i\in\mathcal{C}$, and hence $\mathcal{C}$ covers $U$. Since $\mathcal{C}\subseteq\mathcal{S}$ and $|\mathcal{C}|\leq k$, $\mathcal{C}$ is a set cover of $U$ of size at most $k$.
\end{proof}

Using Proposition~\ref{prop:esi-split}, we can now prove the NP-completeness of \textsc{ESI Set} problem on split graphs.

\begin{theorem}\label{NPsplit}
The \textsc{ESI Set} problem is NP-complete even when restricted to split graphs.
\end{theorem}

\begin{proof}
Given a set $B\subseteq V(G)$ with $|B|\geq \ell$, we can verify in polynomial time whether $B$ is independent and whether every vertex in $N(B)$ has a neighbor in $V(G)\setminus N[B]$. Hence, \textsc{ESI Set} problem belongs to NP.

To prove NP-completeness, we use the reduction from \textsc{Set Cover} described above. The graph $G'$ can be constructed in polynomial time and is a split graph. We set $\ell=m(k+1)-k$.
By Proposition~\ref{prop:esi-split}, $(U,\mathcal{S},k)$ has a set cover of size at most $k$ if and only if $G'$ has an ESI-set of size at least $\ell$. Since \textsc{Set Cover} is NP-complete and \textsc{ESI Set} problem belongs to NP, it follows that \textsc{ESI Set} problem is NP-complete even when restricted to split graphs.
\end{proof}

We end the section with a discussion on when an optimal ESI-set of a split graph can be obtained from an optimal isolating set and vice versa.

The above results show that we can easily construct a split graph $G$ with $\es(G)=k$ for some $k \leq \alpha(G)-1$. On the other hand, $\iota(G)=1$ holds for every split graph $G=K+I$ with at least one edge, since for any $x \in K$ the set $\{x\}$ is an isolating set of $G$. Recall that for any ESI-set $B$ of $G$, the set $V(G)\setminus  N[B]$ is an isolating set of $G$. Conversely, if $D$ is an isolating set of $G$, then $V(G) \setminus  N[D]$ is an ESI-set of $G$. However, an optimal ESI-set does not necessarily correspond to an optimal isolating set, and vice versa. Therefore, we are interested in determining for which split graphs $G$ and which $\es(G)$-sets $B$ the set $V(G)\setminus N[B]$ is an $\iota(G)$-set. Since $V(G) \setminus N[B]$ is always an isolating set, it suffices to determine when $|V(G)\setminus N[B]| = 1,$
as this implies that the resulting isolating set is optimal.

\begin{theorem}\label{thm:ESIvsIOTA}
Let $G=K+I$ be a split graph not isomorphic to a complete graph or an edgeless graph and let $|I|=\alpha(G)$. Then $G$ contains an $\es(G)$-set $B$ such that $V(G) \setminus N[B]$ is an $\iota(G)$-set if and only if $\omega(G)=|K|+1$ and $\min{\{{\rm{deg}}_I(x):\, x \in K\}} \geq 2$.
\end{theorem}

\begin{proof}
First let $B$ be an $\es(G)$-set of a split graph $G=K+I$ such that $V(G)\setminus N[B]$ is an $\iota(G)$-set. By Lemma~\ref{l:insideI} $B \subseteq I$.  Since $|B|=\es(G) \leq \alpha(G)-1$, there exists $x \in I$ such that $x \notin B$. Since $\iota(G)=1$, it follows that $|V(G) \setminus N[B]|=1$ and as $x \notin N[B]$ we deduce that $V(G) \setminus N[B]=\{x\}$, $K \subseteq N(B)$ and $B=I \setminus \{x\}$. By the ES-condition for $B$, it follows that every $y \in N(B)=K$ has a neighbor in $V(G) \setminus N[B]=\{x\}$ which implies that $x$ is adjacent to all vertices of $K$. Thus $\omega(G)=|K|+1$. It is also clear that $\min{\{{\rm{deg}}_I(x):\, x \in K\}} \geq 2$ as every vertex of $K$ has a neighbor in $B$ and it is also adjacent to $x$.

Conversely, let $G$ be a split graph with $\omega(G)=|K|+1$ and $\min{\{{\rm{deg}}_I(z):\, z \in K\}} \geq 2$. Then there exists $x \in I$ such that $x$ is adjacent to all vertices of $K$. Let $B=I \setminus \{x\}$. Since $\min{\{{\rm{deg}}_I(z):\, z \in K\}} \geq 2$, it follows that every vertex $y \in K$ is adjacent to at least one vertex in $B$ and thus $N(B)=K$. Since $N(I \setminus B)=N(\{x\})=K$, Lemma~\ref{l:type2ESIsets} implies that $B$ is an ESI-set of $G$. Since $|B|=\alpha(G)-1$, $B$ is an $\es(G)$-set by Theorem~\ref{alpha_up_bound}.  Moreover, $B$ has the property that $V(G)\setminus N[B]=\{x\}$ and thus $V(G) \setminus N[B]$ is an optimal isolating set.

\end{proof}

To see  for which $\iota(G)$-set $D$, the set $V(G) \setminus N[D]$ is an $\es(G)$-set we first need the following.

\begin{lemma}\label{l:pomozna_split}
Let $G=K+I$ be a split graph not isomorphic to a complete graph with $\alpha(G)=|I|$. If there exists $x \in I$ such that the set $\{x\}$ is an isolating set of $G$, then $N(x)=K$.  
\end{lemma}

\begin{proof}
Let $x \in I$ such that $\{x\}$ is an isolating set of $G$. Suppose that there exists $y$ in $K$ such that $xy \notin E(G)$. Since $\{x\}$ is an isolating set, $G\setminus N[x]$ has no edges and thus $y$ has no neighbors in $I$. This leads to a contradiction with $|I|=\alpha(G)$. 
\end{proof}

\begin{theorem}\label{thm:ESIvsIOTA2}
Let $G=K+I$ be a split graph with $|I|=\alpha(G)$ that is not isomorphic to a complete graph or an edgeless graph. Then there exists an $\iota(G)$-set $D$ such that $V(G) \setminus N[D]$ is an $\es(G)$-set if and only if $\omega(G)=|K|+1$ and $\min{\{{\rm{deg}}_I(x):\, x \in K\}} \geq 2$ or there exists $x \in K$ such that $I \setminus N_I(x)$ is an $\es(G)$-set. 
\end{theorem}
\begin{proof}
Let $D=\{x\}$ be an $\iota(G)$-set and let first $B=V(G) \setminus N[D]$ be an $\es(G)$-set. Then $x \in K$ or $x \in I$ and $N(x)=K$ by Lemma~\ref{l:pomozna_split}. Assume first that $x \in I$ and $N(x)=K$. Then $\omega(G)=|K|+1$. If there exists $y \in K$ with ${\rm{deg}}_I(y)=1$, i.e.\ with $x$ being the only neighbor of $y$ in $I$, then $B=I \setminus N_I(y)$ is an $\es(G)$-set. Otherwise $\omega(G)=|K|+1$ and $\min{\{{\rm{deg}}_I(x):\, x \in K\}} \geq 2$.
 
Assume now that $x \in K$. Then $B = V(G)\setminus N[x]=I \setminus N_I(x) \subseteq I$. Since $|B|=\es(G)$, Lemma~\ref{l:type1ESIsets} implies that ${\rm{deg}}_I(x)=\min{\{{\rm{deg}}_I(x):\, x \in K\}}$ and thus there exists $x \in K$ such that $I\setminus N_I(x)$ is an $\es(G)$-set.

For the converse, let $\omega(G)=|K|+1$ and $\min{\{{\rm{deg}}_I(x):\, x \in K\}} \geq 2$ or there exists $x \in K$ such that $I \setminus N_I(x)$ is an $\es(G)$-set. First assume that  $\omega(G)=|K|+1$ and $\min{\{{\rm{deg}}_I(x):\, x \in K\}} \geq 2$. Then by Theorem~\ref{thm:ESIvsIOTA} there exists an $\es(G)$-set $B$ such that $V(G) \setminus N[B]$ is $\iota(G)$-set. Hence Proposition~\ref{optimalToOptimal} implies that there exists an $\iota(G)$-set $D$ such that $V(G) \setminus N[D]$ is an $\es(G)$-set.

Finally, assume that there exists $x \in K$ such that $B=I \setminus N_I(x)$ is an $\es(G)$-set. Then $\{x\}$ is an $\iota(G)$-set with the property that $B=V(G) \setminus N[x]$ is an $\es(G)$-set.
 \end{proof}

\section{Proper interval graphs}\label{s:interval}

In this section we present an algorithm that generates a maximum externally supported independent set of an arbitrary proper interval graph.

Recall that a proper interval graph is an interval graph that has an intersection model in which no interval properly contains another. Proper interval graphs are exactly unit interval graphs~\cite{roberts1969, gardi2007}, which means that they have an interval representation in which every interval has unit length. Equivalently, proper interval graphs are exactly $K_{1,3}$-free interval graphs~\cite{roberts1969}. We may assume without loss of generality that all endpoints of the unit interval representation are distinct \cite{gt2004}. Moreover, we may assume that proper interval graphs considered in this section are connected. If not, we can use the below presented algorithm separately on each component.

Before describing the algorithm, we preprocess the unit interval
representation of a proper interval graph $G$ with $v_i=[a_i,b_i]$, $i=1,\ldots,n$, such that $a_1 < a_2 < \ldots < a_n$.

For every $k\in\{1,\ldots,n\}$, let
$$p_a(k)=\max\{j:b_j<a_k\},$$
$$p_b(k)=\max\{j:a_j<b_k\},$$
$$s_b(k)=\min\{j:a_j>b_k\}.$$

We set $p_a(k)=0$ and $s_b(k)=n+1$ whenever the corresponding sets are empty.
The following preprocessing algorithm computes all values $p_a(k)$, $p_b(k)$, and $s_b(k)$ in $O(n)$ time using two-pointer scans. After preprocessing, each of
these values can be accessed in constant time. 

\medskip

\begin{algorithm}[H]
\caption{Preprocessing of the unit interval representation}
\label{alg:preprocessing}
\KwIn{A proper interval graph $G$ with a unit interval representation
      $v_i=[a_i,b_i]$, $i=1,\ldots,n$, with $a_1< \ldots < a_n$}.
\KwOut{The arrays $p_a$, $p_b$, and $s_b$.}

$j \gets 0$\;
\For{$k=1$ \KwTo $n$}{
    \While{$j<n$ and $b_{j+1}<a_k$}{
        $j\gets j+1$\;
    }
    $p_a(k)\gets j$\;
}

$j \gets 0$\;
\For{$k=1$ \KwTo $n$}{
    \While{$j<n$ and $a_{j+1}<b_k$}{
        $j\gets j+1$\;
    }
    $p_b(k)\gets j$\;
}

$j \gets 1$\;
\For{$k=1$ \KwTo $n$}{
    \While{$j\leq n$ and $a_j\leq b_k$}{
        $j\gets j+1$\;
    }
    $s_b(k)\gets j$\;
}
\Return{$p_a,p_b,s_b$}\;
\end{algorithm}

\medskip 
We now use the precomputed values $p_a$, $p_b$, and $s_b$ to find a maximum ESI-set.
The algorithm is given below.

\medskip

\begin{algorithm}[H]
\caption{Maximum ESI-set of a proper interval graph}
\label{alg:esi-proper}
\KwIn{A proper interval graph $G$ with a unit interval representation
      $v_i=[a_i,b_i]$, $i=1,\ldots,n$, ordered from left to right.}
\KwOut{A maximum ESI-set $B$.}

$B \gets \emptyset$\;
$k \gets 1$\;

\While{$k \leq n$}{
    \If{
    $(k=1$ \Or ($p_a(k)\geq 1$ \And $a_{k-1}\leq b_{p_a(k)}))$ \And\\
    \hspace{1em} $(k=n$ \Or ($s_b(k)\leq n$ \And $a_{s_b(k)}\leq b_{k+1}))$
    }   {
        $B \gets B\cup\{v_k\}$\;
        $i \gets p_b(k)$\;

        \If{$i=n$ \Or $a_n < b_{i+1}$}{
            \Return{$B$}\;
        }
        \Else{
            $k \gets s_b(i+1)$\;
        }
    }
    \Else{
        $k \gets k+1$\;
    }
}
\Return{$B$}\;
\end{algorithm}


\medskip

Before proving the correctness of the algorithm, we prove the following lemma.

\begin{lemma} \label{lemma:algoritem}
Let $G$ be a proper interval graph and let $B$ be an ESI-set of $G$.
Then no two distinct vertices of $B$ have a common neighbor.
\end{lemma}

\begin{proof}
Suppose that two distinct vertices $u,v\in B$ have a common neighbor $w$. By the ES-condition, there exists $z\in V(G)\setminus N[B]$ adjacent to $w$. Since $B$ is independent and $z\notin N[B]$, the vertices $u$, $v$, and $z$ are pairwise nonadjacent. Thus, the subgraph of $G$ induced by vertices $u$, $v$, $w$, and $z$ is isomorphic to $K_{1,3}$, a contradiction since proper interval graphs are $K_{1,3}$-free.
\end{proof}

\begin{theorem}
For a proper interval graph $G$, Algorithm~\ref{alg:esi-proper} returns an $\es(G)$-set $B$ in $O(n)$ time and $O(n)$ space.
\end{theorem}

\begin{proof}
Let $G$ be a proper interval graph with unit interval representation $v_i=[a_i,b_i]$ for $i \in [n]$ and $a_1< a_2 < \ldots < a_n$. 

First, we show that the set $B$ constructed with Algorithm~\ref{alg:esi-proper} is an ESI-set. The set $B$ is clearly independent, since after adding a vertex $v_k$ to $B$, the algorithm skips all neighbors of $v_k$.
It remains to verify the ES-condition. Let $v_k\in B$ and let $x\in N(v_k)$. If $x$ lies to the left of $v_k$, then $v_{p_a(k)}$ is adjacent to $x$, while if $x$ lies to the right of $v_k$, then $v_{s_b(k)}$ is adjacent to $x$. By the definitions of $p_a(k)$ and $s_b(k)$, both vertices $v_{p_a(k)}$ and $v_{s_b(k)}$ lie outside $N[v_k]$. 
Once $v_k$ is added to $B$, the algorithm skips its neighbors and the neighbors of $v_{s_b(k)}$ which is the first vertex to its right outside the closed neighborhood of $v_k$. Consequently, neither $v_{p_a(k)}$ nor $v_{s_b(k)}$ is adjacent to any other vertex of $B$. Therefore, every $x\in N(v_k)$ has a neighbor in $V(G)\setminus N[B]$. This shows that $B$ satisfies the ES-condition.

It remains to prove that $|B|=\es(G)$. Let $D$ be an $\es(G)$-set such that $|B\cap D|$ is as large as possible. Suppose to the contrary that $B\neq D$. The algorithm constructs the set $B$ greedily from left to right until there is no possible vertex to add into the set $B$. Therefore, the set $B$ is maximal and it is not possible that $B \subset D$. Now, let $v_k$ be the first vertex, from left to right, where $B$ and $D$ differ. Since $B$ and $D$ coincide on all vertices preceding $v_k$, and the algorithm always selects the first possible vertex from left to right, we have $v_k\in B\setminus D$. There must be a vertex of $D\setminus B$ to the right of $v_k$, otherwise $D\subset B$, implying $|B|>|D|=\es(G)$, a contradiction. Let $v_{k_1}$ be the first (leftmost) such vertex. 
Denote by $D_{\ell} = \{v_i : v_i \in D \land a_i < a_{k_1}\}$  and $D_r = \{v_i : v_i \in D \land a_i > a_{k_1}\}$ the vertices of $D$ to the left and to the right of $v_{k_1}$, respectively. By Lemma~\ref{lemma:algoritem}, the sets $N[D_l]$, $N[v_{k_1}]$ and $N[D_r]$ are pairwise disjoint.
Define $D'=(D\setminus \{ v_{k_1}\})\cup \{v_k\} = D_{\ell} \cup  \{ v_{k}\} \cup D_r $. We show that $D'$ is an ESI-set. 
Since all vertices of $D_\ell$ are also in the independent set $B$ and as $v_k$ is also in $B$, no vertex of $D_\ell$ is adjacent to $v_k$. By the ordering of the intervals, no vertex of $D_r$ is adjacent to $v_k$ either. Thus, $D'$ is independent. It remains to verify the ES-condition. 

If $D_r = \emptyset$, then $D'\subseteq B$, and hence $D'$ is an ESI-set. Thus, suppose $D_r \neq \emptyset$.
To the left of $v_k$, the ES-condition holds by the choice of $D'$. To the right of $v_k$, the vertex $v_{s_b(k)}$ is adjacent to every
right neighbor of $v_k$. We show that $v_{s_b(k)}\notin N[D_r]$. 
Suppose first that $k_1<s_b(k)$. Since the algorithm selects $v_k$, we have $a_{s_b(k)}\leq b_{k+1}\leq b_{k_1}$. Thus, $v_{s_b(k)}$ is adjacent to $v_{k_1}$. By Lemma~\ref{lemma:algoritem}, it is not adjacent to any vertex of $D_r$. Suppose now that $s_b(k)\leq k_1$. Since $D$ is independent, no vertex of $D_r$ is adjacent to $v_{k_1}$. By the ordering of the intervals, no vertex of $D_r$ is adjacent
to $v_{s_b(k)}$ either. Therefore, $v_{s_b(k)}\notin N[D_r]$.
The ES-condition remains satisfied for the vertices of $D_r$.
Indeed, any external support of a vertex in $N(D_r)$
that is adjacent to $v_k$ would also be adjacent to
$v_{k_1}$, by the ordering of the intervals.
This contradicts the fact that $D$ is an ESI-set.


Since $|D'|=|D|=\es(G)$, the set $D'$ is an $\es(G)$-set. Moreover, since $v_{k_1} \notin B$ and $v_k\in B$, we have $|B\cap D'|>|B\cap D|$, which contradicts the choice of $D$ as an $\es(G)$-set with $|B\cap D|$ as large as possible. Therefore, $|B|=\es(G)$.

It remains to show that both the time and space complexity are $O(n)$.
In Algorithm~\ref{alg:esi-proper}, the index $k$ only increases,
and each iteration takes constant time. Thus, the algorithm runs
in $O(n)$ time. Together with the preprocessing, the total running time is $O(n)$.
The three preprocessing arrays require $O(n)$ space, and the set $B$
contains at most $n$ vertices. Thus, the space complexity
is also $O(n)$.
\end{proof}


\section{Conclusion}

Our results exhibit two different computational behaviors within structurally restricted graph classes of unbounded clique-width. On split graphs, the \textsc{ESI Set} problem is computationally hard, whereas on proper interval graphs it admits a linear-time algorithm. This contrast motivates the study of the full class of interval graphs. In particular, it is natural to ask whether the linear-time approach for proper interval graphs can be extended to arbitrary interval graphs, or whether the increased flexibility of general interval representations leads to a different computational complexity. This gives rise to the following question:

\begin{question}
Is the \textsc{ESI Set} problem solvable in polynomial time on interval graphs?
\end{question}

\vskip 1pc \noindent{\bf Acknowledgments.} 
This work has been supported by the European Commission's Horizon Europe Research and Innovation programme through the Marie Skłodowska-Curie Actions Staff Exchanges (MSCA-SE) under Grant Agreement no.101182819 (COVER: (C)ombinatorial (O)ptimization for (V)ersatile Applications to (E)merging u(R)ban Problems).
T.D. was partially supported by the Slovenian Research and Innovation Agency (ARIS), program no.~P1--0297 and project number N1-0431. 
D.B. was partially supported by the Slovenian Research and Innovation Agency (ARIS), program no.~P2--0065 and project no.~J1--70016.
A.K. was partially supported by the Slovenian Research and Innovation Agency (ARIS), program no.~P1--0297 and project no.~J1--70016.\\

\vskip 1pc \noindent{\bf Statements and Declarations.}
Authors have no conflict of interest to declare.

\bibliographystyle{plain}
\bibliography{references}

\end{document}